\documentclass[11pt]{amsart}

\newcommand{\R}{\mathbb R}

\newcommand{\Q}{\mathbb Q}
\newcommand{\Z}{\mathbb Z}

\newcommand{\p}{\mathbf p}
\newcommand{\q}{\mathbf q}
\renewcommand{\r}{\mathbf r}

\theoremstyle{plain}
\newtheorem{Thm}{Theorem}
\newtheorem{Prop}[Thm]{Proposition}

\theoremstyle{definition}

\newtheorem{Rmk}{Remark}
\newtheorem{Exam}{Example}

\begin{document}

\title{On the mass center of the tent map}
\author[Kuo-Chang Chen and Xun Dong]{Kuo-Chang Chen and Xun Dong}
\address{Kuo-Chang Chen, Department of Mathematics \\National Tsing Hua University \\Hsinchu\\Taiwan}
\email{kchen@math.nthu.edu.tw}
\address{Xun Dong, Department of Mathematics \\University of Miami \\Coral Gables, FL 33124\\ USA}
\email{xundong@math.miami.edu}
\date{February 18.}

\begin{abstract}
It is well known that the time average or the center of mass for
generic orbits of the standard tent map is $0.5$. In this paper we
show some interesting properties of the exceptional orbits,
including periodic orbits, orbits without mass center, and orbits
with mass centers different from $0.5$. We prove that for any
positive integer $n$, there exist $n$ distinct periodic orbits for
the standard tent map with the same center of mass, and the set of
mass centers of periodic orbits is a dense subset of $[0,2/3]$.
Considering all possible orbits, then the set of mass centers is
the interval $[0,2/3]$. Moreover, for every $x$ in $[0,2/3]$,
there are uncountably many orbits with mass center $x$. We also
show that there are uncountably many orbits without mass center.
\end{abstract}
\maketitle

\section{Introduction}\label{sec:intr}

Let $T:[0,1]\rightarrow[0,1]$ be the standard tent map defined by
$T(x)=1-|2x-1|$ or
$$
  T(x)=\left\{\begin{array}{ll}
               2x & \textrm{when}\; x\in[0,\frac{1}{2}]\smallskip\\
               2(1-x) & \textrm{when}\; x\in[\frac{1}{2},1].
               \end{array}\right.
$$
This map
is often introduced as the first example of chaotic maps in typical textbooks for
dynamical systems.
Its dynamics exhibit various features that are commonly used to identify chaotic systems.
Some of these well-known features are sensitive dependence,
topological transitivity (for almost all orbits in the Lebesgue sense),
existence of infinitely many orbits with positive Lyapunov exponent,
and being a topological factor of every unimodal map with topological entropy $\log 2$.
Details and further references can be found in \cite{KH}, for instance.

The $n$-th iterate of $T$ is given by
$$
  T^n(x)=\left\{\begin{array}{ll}
                 2^n\big(x-\frac{2m}{2^n}\big) & \textrm{when}\; x\in[\frac{2m}{2^n},\frac{2m+1}{2^n}]\smallskip\\
                 2^n\big(\frac{2(m+1)}{2^n}-x\big) & \textrm{when}\; x\in[\frac{2m+1}{2^n},\frac{2(m+1)}{2^n}]
                 \end{array}\right.
$$
for $m\in\{0,1,\cdots,2^{n-1}-1\}$.
Fixed points of $T^n$, or $n$-periodic points of $T$, are clearly of the form
$\frac{2m}{2^n-1}$ or $\frac{2m}{2^n+1}$.
Given an $n$-periodic point $x\in(0,1)$,
the average $$\bar{x}=\frac{1}{n}\left(x+T(x)+\cdots+T^{n-1}(x)\right)$$ over its orbit
is called the {\em center of mass} of the orbit.

In \cite{Mis} Misiurewicz noted that two cycles with mirror
itineraries have the same center of mass. He then raised an
interesting question: Can there be three different cycles with the
same center of mass? The answer is affirmative as can be easily
verified from the distinct orbits of
$$
  \frac{166}{4095},\;\frac{202}{4095},\;\frac{278}{4095},\;\frac{418}{4095};
$$
or the distinct orbits of
$$
  \frac{42}{4097},\;\frac{54}{4097},\;\frac{134}{4097},\;\frac{190}{4097}.
$$
The orbits of
$$
  \frac{2250}{16385},\;\frac{2266}{16385},\;\frac{2446}{16385},\;\frac{2490}{16385},\;
  \frac{2510}{16385}
$$
are also distinct but with the same mass center, and so are the orbits of
$$
  \frac{1510}{65535},\;\frac{1658}{65535},\;\frac{2270}{65535},\;\frac{3566}{65535},\;
  \frac{3830}{65535},\;\frac{3938}{65535}.
$$
As the denominator increases there are more and more examples of distinct cycles with the
same mass center.
The search for $n$ distinct cycles with the same mass center becomes a challenging task
when $n$ is large.
In this paper we prove that

\begin{Thm}\label{main theorem}
For any positive integer $n$, there exist $n$ distinct cycles for
the standard tent map with the same center of mass.
\end{Thm}

The definition for center of mass can be extended to non-periodic orbits.
For any $x\in[0,1]$, if the partial time average $\frac{1}{n}\sum_{k=0}^{n-1}T^ix$ converges as
$n\rightarrow\infty$, the limit is called the {\em center of mass}
or the {\em time average}
of the orbit.
This definition clearly coincides with the definition on periodic orbits.
Since the Lebesgue measure is an ergodic invariant measure,
by the Birkhoff ergodic theorem, we know that for almost every $x$
(in the sense of Lebesgue) the mass center of its orbit is $\frac{1}{2}$, which is
the space average of $T$.
Also, for almost every $x$ the frequency of appearance of its forward orbit in
a measurable set $A$ is equal to the Lebesgue measure of $A$.

For interval maps with an absolutely continuous ergodic invariant
measure, in literature there is little study on the collection of
exceptional orbits, such as orbits without mass center, and orbits
with mass centers different from the space average. These orbits
are also of our interest because they exhibit dynamics of the map
with different degrees of randomness. What can we say about the
set of orbits without mass center, except for being a zero measure
set? We will show that

\begin{Thm}\label{thm:nocenter}
There are uncountably many orbits of the standard tent map without mass center.
\end{Thm}

An interesting comparison with the tent map is the irrational rotation which maps $x\in[0,1]$
to $x+\alpha\;\,\textrm{mod}\; 1$ for some fixed $\alpha\in\R\setminus\Q$.
Irrational rotations also has the Lebesgue measure as an ergodic invariant measure
but every orbit is equidistributed and thus has $\frac{1}{2}$ as its mass center (see \cite{AA});
in other words, there is no exceptional orbit for irrational rotations.

Counting all possible orbits, periodic and non-periodic ones, how
do their mass centers distribute over $[0,1]$? Consider the set of
all possible mass centers. It might appear that this set is a
meagre set, but it is indeed not the case. We will show in this
paper that

\begin{Thm}\label{thm:dense}
The set of all mass centers for the standard tent map is the
interval $[0,\frac{2}{3}]$. The set of mass centers for cycles is
a dense subset of this interval. Moreover, for every $x$ in
$[0,\frac{2}{3}]$, there are uncountably many orbits with mass
center $x$.
\end{Thm}

In section~\ref{sec:binary} we provide a sufficient condition for two periodic orbits
to have the same mass center.
Section~\ref{sec:family} is devoted to the proof of Theorem~\ref{main theorem}.
Theorem~\ref{thm:nocenter} and Theorem~\ref{thm:dense} are both proved
in section~\ref{sec:masscenters}.
Several examples and related problems can be found in section~\ref{sec:ex} and \ref{sec:rmk}.

\section{Binary representations}\label{sec:binary}

In this section we introduce the concept of trace vector and establish
a sufficient condition for periodic orbits of certain type to have
the same mass center.

Throughout this and the next section
we consider $n$-periodic points of the form $x=\frac{2m}{2^n-1}$.
The numerator $2m$ of $x$ can be uniquely written as
$$a_1 2^{n-1}+a_2 2^{n-2}+\cdots+a_{n-1}2^1$$ with $a_i\in\{0,1\}$ for all $i$.
For convenience we add the term $a_n 2^0$ with coefficient $a_n=0$ and
write the numerator as a vector in $(\Z/2\Z)^n$:
$$
2m = [a_1,a_2,\cdots,a_n].
$$
The tent map $T$ acts on the vector as a cyclic left-shift when
$a_1=0$, while $0$'s and $1$'s are swapped after a cyclic
left-shift when $a_1=1$.
More precisely, the tent map induces a map $T_*$ 
from $\{0,1,2,\cdots,2^n-1\}$ to $\{0,2,4,\cdots,2^{n-1}-2\}$ given by
\begin{eqnarray*}
  T_*[a_1,a_2,\cdots,a_n] &=&
  \left\{\begin{array}{ll}
          \left[a_2,a_3,\cdots,a_n,a_1\right] & \textrm{if}\; a_1=0\smallskip\\
          \left[1-a_2,1-a_3,\cdots,1-a_n,1-a_1\right] & \textrm{if}\; a_1=1.
         \end{array}\right. \\
  &=& [a_2,a_3,\cdots,a_n,a_1] + [a_1,a_1,\cdots,a_1,a_1] \\
  &=& [a_1,a_2,\cdots,a_{n-1},a_n] E F,
\end{eqnarray*}
where $E, F \in M_{n \times n} (\Z/2\Z)$:
$$
E = \left[
        \begin{array}{ccccccc}
        0&0&\cdots &0&1\\
        1&0&\cdots &0&0\\
        0&1&0&\cdots &0\\
        \vdots&\ddots&\ddots&\ddots&\vdots\\
        0&\cdots &0&1 &0
        \end{array}
    \right], \;\;
F = I_n +
    \left[
        \begin{array}{ccccccc}
        0&0&\cdots &0&0\\
        0&0&\cdots &0&0 \\
        \vdots&\vdots&\ddots&\vdots&\vdots\\
        0&0&\cdots &0&0\\
        1&1&\cdots &1&1
        \end{array}
    \right].
$$
Note that the last entry of $T_*[a_1,a_2,\cdots,a_n]$ is always zero.
Therefore,
\begin{eqnarray*}
  T_*^k[a_1,a_2,\cdots,a_n]
  &=& [a_1,a_2,\cdots,a_n] (E F)^k \\
  &=& [a_1,a_2,\cdots,a_n] E^k(E^{1-k}FE^{k-1})(E^{2-k}FE^{k-2})\cdots(E^{-1}FE^1)F \\
  &=& [a_{k+1},\cdots,a_n,a_1,\cdots,a_k]+[a_k,a_k,\cdots,a_k].
\end{eqnarray*}
The last identity holds because
$$[a_1,a_2,\cdots,a_n] E^k=[a_{k+1},\cdots,a_n,a_1,\cdots,a_k]$$
and the collected effect of
$$(E^{1-k}FE^{k-1})(E^{2-k}FE^{k-2})\cdots(E^{-1}FE^1)F$$
is adding a vector $[c,c,\cdots,c]$ to $[a_{k+1},\cdots,a_n,a_1,\cdots,a_k]$.
Since the last entry of $T_*^k[a_1,a_2,\cdots,a_n]$ is zero, $c$ must be equal to $a_k$.

The orbit of $x$ is encoded in the $n$ by $n$ matrix $A(x)\in M_{n\times n}(\Z/2\Z)$:
\begin{eqnarray*}
A(x) &=& \left[
           \begin{array}{c}
            a_1,a_2,\cdots,a_n\\
            T_*[a_1,a_2,\cdots,a_n] \\
            T_*^2[a_1,a_2,\cdots,a_n] \\
            \vdots \\
            T_*^{n-1}[a_1,a_2,\cdots,a_n]
           \end{array}
         \right]  \\
     &=& \left[
           \begin{array}{ccccc}
            a_1&a_2&\cdots&a_{n-1}&a_n\\
            a_2&a_3&\cdots&a_n&a_1\\
            \vdots & \vdots & \ddots & \vdots & \vdots \\
            a_n&a_1&\cdots&a_{n-2}&a_{n-1}
           \end{array}
         \right] + \left[
           \begin{array}{cccc}
            a_n&a_n&\cdots&a_n\\
            a_1&a_1&\cdots&a_1\\
            \vdots & \vdots & \ddots & \vdots \\
            a_{n-1}&a_{n-1}&\cdots&a_{n-1}
           \end{array}
         \right]\\
     &=& \left[\alpha,E^{-1}\alpha,\cdots,E^{1-n}\alpha\right]
         +\left[E\alpha,E\alpha,\cdots,E\alpha\right],
\end{eqnarray*}
where $\alpha=[a_1,a_2,\cdots,a_n]^T$. Let $\beta = \alpha + E\alpha$ and
$$
 C_k = I+E^{-1}+\cdots+E^{-k}.
$$
Then
\begin{eqnarray}\label{eqn:ab}
  E^{-k}\alpha+E\alpha = (I+E^{-1}+\cdots+E^{-k})(\alpha + E\alpha) = C_k\beta.
\end{eqnarray}
Therefore
\begin{eqnarray*}\label{eqn:A}
  A(x)=\left[\beta,C_1\beta,\cdots,C_{n-1}\beta\right].
\end{eqnarray*}
The last column of $A(x)$ is zero. Since all entries of the matrix
$C_{n-1}$ are $1$, we conclude that the number of $1$'s in $\beta$
must be even. Also note that {\em the first column
$\beta=\alpha+E\alpha$ of $A(x)$ is exactly the itinerary of
$x=\frac{2m}{2^n-1}$}; that is,
$\beta_i=\chi_{[\frac{1}{2},1]}(T^{i-1}(x))$ for each $i$.

Note that $\frac{2m}{2^n-1}$ and $\frac{2k}{2^n-1}$ can not have
the same itinerary unless $m=k$. Since there are $2^{n-1}$ vectors
of the form $[a_1,\cdots,a_{n-1},0]$ and $2^{n-1}$ possible
itineraries $\beta$ (which must have evenly many $1$'s), the
correspondence between $x=\frac{2m}{2^n-1}$ and $\beta$ is
necessarily bijective.

Let $s_j(\beta)=s_j$ be the number of $1$'s in the $j$-th column
$C_{j-1} \beta$ of $A(x)$. We call
$$
  s(\beta) = (s_1,s_2,\cdots,s_n)\in\Z^n
$$
the {\em trace vector} of $x$. If follows easily from
(\ref{eqn:ab}) that $s_n=0$ and each $s_j$ is even.
Since
\begin{eqnarray*}
  C_{n-k-1} \beta =E^{1-(n-k)}\alpha+E\alpha = E^{1+k}\alpha+E\alpha
  = E^k(E\alpha+E^{1-k}\alpha) = E^k C_{k-1} \beta
\end{eqnarray*}
we also have
\begin{eqnarray}\label{eqn:sym}
s_{n-k} = s_k\;\;\textrm{for}\;\; k=1,\cdots,n-1.
\end{eqnarray}

Write $A(x)=[b_{ij}]$ with indices $i$ and $j$ ranging from $1$ to $n$,
then $b_{ij}=\sum_{r=0}^{j-1} \beta_{i+r}$. It is understood that the
subindex $i+r$ of each $\beta_{i+r}$ is in $\Z/n\Z$.
The orbit of $x$ has mass center
\begin{eqnarray}\label{eqn:cm}
  \bar{x} &=& \frac{1}{(2^n-1)n}\sum_{i,j=1}^{n} 2^{n-j} b_{ij}
       \;\;=\;\;\frac{1}{(2^n-1)n}\sum_{j=1}^{n} 2^{n-j} s_j.
\end{eqnarray}
Thus {\em for two $n$-periodic points $x=\frac{2m}{2^n-1}$ and
$y=\frac{2k}{2^n-1}$, a sufficient condition for $\bar{x}=\bar{y}$
is that they have the same trace vector}.

Let $\beta$ and $\gamma$ denote respectively the itineraries of
$n$-periodic points $x$ and $y$. As explained earlier, $x=y$ if
and only if $\beta=\gamma$. It follows that $x$ and $y$ are in the
same cycle if and only if $\gamma = E^k \beta$ for some $k$, in
which case we will say that $\beta$ and $\gamma$ are {\em
equivalent}.

Let
$$
D = \left[
        \begin{array}{ccccccc}
        0&\cdots&0&0&1\\
        0&\cdots&0&1&0\\
        \vdots& && &\vdots\\
        0&1&0&\cdots &0\\
        1&0&0&\cdots &0
        \end{array}
    \right].
$$
Then $D$ and $E$ generate a dihedral group of order $2n$ with
relations
\begin{eqnarray}\label{eqn:dihedral}
  E^n=D^2=(ED)^2=I_n.
\end{eqnarray}
Observe that
$$
C_kE=EC_k, \;\; C_kD=DE^kC_k
$$
for all $k$.

If $\gamma=E\beta$, then $C_k\gamma=C_k E\beta=E C_k\beta$ for all $k$.
Thus $A(x)$ and $A(y)$ have the same trace vector as expected.

If $\gamma=D\beta$, then $C_k\gamma=C_k D\beta=D E^k C_k \beta$
for all $k$. Thus $A(x)$ and $A(y)$ have the same trace vector.
Note that $\gamma$ is the mirror image of $\beta$. In this case
$x$ and $y$ belong to different cycles except when
$\beta=DE^k\beta$ for some $k$. For instance, the point
$x=\frac{26}{127}$ has itinerary
$$
 \beta = [0,0,1,0,1,1,1]^T.
$$
The itinerary of $y=\frac{88}{127}$ is $\beta$ reversed.
These two points $x$ and $y$ belong to two different cycles with the same center of mass $\frac{72}{127}$.

One might attempt to look for other permutation matrices $P$ with the nice property
that $s(\beta) = s(P\beta)$ for all itineraries $\beta$.
Unfortunately there are no such permutation matrices other than those in
the dihedral group generated by $D$ and $E$.
In order to find three or more different cycles with the same mass center,
we need more explicit constructions, as to be shown in the next section.

\section{A family of cycles with the same mass center}\label{sec:family}

In order to construct $n$ distinct cycles with the same mass
center, we will find $n$ $(8n)$-periodic points such that their
itineraries $\beta^{(1)},\cdots\beta^{(n)}$ have the same trace
vector, yet $\beta^{(i)}$ and $\beta^{(j)}$ are non-equivalent
when $i \neq j$. The ``$n$'' in this section is not to be confused
with the ``$n$'' in the previous section, which is replaced by
``$8n$'' here for convenience.

For each $1\leq m\leq n$, let
\begin{eqnarray} \label{eqn:betam}
 \beta^{(m)} & = & [\beta^{(m)}_1,\beta^{(m)}_2,\cdots,\beta^{(m)}_{8n}]^T
\end{eqnarray}
be defined by the following conditions
\begin{eqnarray*}
 \beta^{(m)}_{i}+\beta^{(m)}_{4n+i} &=& 1 \quad\textrm{for}\; i=1,2,\cdots,4n, \\
 \beta^{(m)}_{i} &=& \left\{
                        \begin{array}{ll}
                          1 & \;\textrm{if}\;\; i\in\{2m,2m+1\}\\
                          0 & \;\textrm{if}\;\; i\in\{1,2,\cdots,4n\}\setminus\{2m,2m+1\}.
                        \end{array}\right.
\end{eqnarray*}
It is easy to see that these conditions uniquely determine
$\beta^{(m)}$. It is understood that the subindex $i$ of each
$\beta^{(m)}_{i}$ is in $\Z/8n\Z$. For instance, when $n=3$, the
$\beta^{(m)}$'s are
\begin{eqnarray}
 \beta^{(1)} & = & [0,1\,|\,1,0\,|\,0,0\,|\,0,0\,|\,0,0\,|\,0,0\,|\,1,0\,|\,0,1\,|\,1,1\,|\,1,1\,|\,1,1\,|\,1,1]^T  \nonumber\\
 \beta^{(2)} & = & [0,0\,|\,0,1\,|\,1,0\,|\,0,0\,|\,0,0\,|\,0,0\,|\,1,1\,|\,1,0\,|\,0,1\,|\,1,1\,|\,1,1\,|\,1,1]^T \label{eqn:beta123}\\
 \beta^{(3)} & = & [0,0\,|\,0,0\,|\,0,1\,|\,1,0\,|\,0,0\,|\,0,0\,|\,1,1\,|\,1,1\,|\,1,0\,|\,0,1\,|\,1,1\,|\,1,1]^T  \nonumber
\end{eqnarray}
The vertical bars are inserted to the vectors above to improve readability.

By (\ref{eqn:sym}), we only need to show that $s_k(\beta^{(m)})$
is independent of $m$ for $1\leq k \leq 4n$ because of the
symmetry of the trace vector.

If $k$ is odd, then for each $i\in \Z/8n\Z$ we have
\begin{eqnarray*}
  b_{i,k}^{(m)}+b_{4n+i,k}^{(m)}=
  \sum_{j=0}^{k-1} \beta^{(m)}_{i+j} + \sum_{j=0}^{k-1} \beta^{(m)}_{4n+i+j}
  = k \equiv 1 \mod 2.
\end{eqnarray*}
Hence there are as many 1's as 0's in $[b_{1,k}^{(m)}\,
b_{2,k}^{(m)} \cdots b_{8n,k}^{(m)}]^T$, and therefore
$s_k(\beta^{(m)})=4n$ when $k$ is odd.

Next we note that
$$
b_{i,4n}^{(m)}=\sum_{j=0}^{4n-1} \beta_{i+j}^{(m)}=
b_{i+1,4n}^{(m)} +\beta_i^{(m)} -\beta_{4n+i}^{(m)} = b_{i+1,
4n}^{(m)} +1.
$$
Hence there are as many 1's as 0's in $[b_{1,4n}^{(m)}\,
b_{2,4n}^{(m)} \cdots b_{8n,4n}^{(m)}]^T$, and therefore
$s_{4n}(\beta^{(m)})=4n$.

We now assume that $k=2r$ where $1\leq r <2n$. Note that
\begin{eqnarray*}
  \beta^{(m)}_{2i-1}+\beta^{(m)}_{2i} &=&
     \left\{\begin{array}{ll}1&\;\textrm{when}\;\;i\in\{m,m+1,2n+m,2n+m+1\},\\
                             0&\;\textrm{otherwise}, \end{array}
     \right. \\
  \beta^{(m)}_{2i}+\beta^{(m)}_{2i+1} &=&
     \left\{\begin{array}{ll}1&\;\textrm{when}\;\;i\in\{0,2n\},\\
                             0&\;\textrm{otherwise}. \end{array}
     \right. \\
\end{eqnarray*}
It follows that
\begin{eqnarray*}
   b^{(m)}_{2i-1,2r}
     &=&\sum_{j=0}^{r-1}\left(\beta^{(m)}_{2i+2j-1}+\beta^{(m)}_{2i+2j}\right)\\
     &=&
     \left\{\begin{array}{ll}1&\;\textrm{if}\;\;
                               i\in\{m+1,m-r+1,2n+m+1,2n+m-r+1\},\\
                             0&\;\textrm{otherwise}, \end{array}
     \right.
\end{eqnarray*}
and
\begin{eqnarray*}
  b^{(m)}_{2i,2r}
     &=&\sum_{j=0}^{r-1}\left(\beta^{(m)}_{2i+2j}+\beta^{(m)}_{2i+2j+1}\right) \\
     &=&
     \left\{\begin{array}{ll}1&\;\textrm{if either}\;
                               \{0,1\}\;\textrm{or}\;\{4n,4n+1\}\subset\{2i,2i+1,\cdots,2i+2r-1\},\\
                             0&\;\textrm{otherwise}. \end{array}
     \right.
\end{eqnarray*}
There are four 1's in the first case and $2r$ 1's in the second
case. Hence $s_k(\beta^{(m)}) = 2r + 4 = k+4$.

It remains to show that if $i \neq j$ then $\beta^{(i)}$ and
$\beta^{(j)}$ are not equivalent; that is, $\beta^{(i)}$ can not
be obtained from $\beta^{(j)}$ by a cyclic rotation. This follows
from the fact that if we write $\beta^{(m)}$ around a circle, then
the length of the longest substring of consecutive 0's is
$4n-2m-1$. Thus the $\beta^{(m)}$'s are itineraries of points in
$n$ distinct cycles with the same center of mass. The proof of
Theorem~\ref{main theorem} is now concluded.

\begin{Rmk}
It is easy to see that the definition of $\beta^{(m)}$ can be
extended to all $1\leq m < 2n$. In fact, all these $2n-1$
itineraries are pairwise non-equivalent. Moreover there is another
itinerary with the same trace vector, making a total of $2n$
distinct $(8n)$-cycles with the same mass center.

Replacing $2n$ by $n$, we may construct $n$ distinct $(4n)$-cycles
with the same mass center as follows. For each $0\leq m < n$, let
\begin{eqnarray} \label{eqn:gammam}
 \gamma^{(m)} & = & [\gamma^{(m)}_1,\gamma^{(m)}_2,\cdots,\gamma^{(m)}_{4n}]^T
 \in(\Z/2\Z)^{4n}
\end{eqnarray}
where $ \gamma^{(m)}_{i}+\gamma^{(m)}_{2n+i} = 1$ for
$i=1,2,\cdots,2n$. For $m=0$ define
\begin{eqnarray*}
 \gamma^{(0)}_{i} &=& \left\{
                        \begin{array}{ll}
                          1 & \;\;\textrm{if}\; i\in\{1,3\}\\
                          0 & \;\;\textrm{if}\; i\in\{1,2,\cdots,2n\}\setminus\{1,3\}
                        \end{array}\right.
\end{eqnarray*}
and for $1\leq m <n$ define
\begin{eqnarray*}
 \gamma^{(m)}_{i} &=& \left\{
                        \begin{array}{ll}
                          1 & \;\;\textrm{if}\; i\in\{2m,2m+1\}\\
                          0 & \;\;\textrm{if}\; i\in\{1,2,\cdots,2n\}\setminus\{2m,2m+1\}.
                        \end{array}\right.
\end{eqnarray*}
Then a similar calculation shows that $s(\gamma^{(m)})$ is
independent of $m$, with $s_{2r-1}=s_{2n}=2n$ and $s_{2r}= 2r+4$
for $r=1,\cdots,n$. It is also not hard to verify that if $i\neq
j$ then $\gamma^{(i)}$ and $\gamma^{(j)}$ are not equivalent.
Therefore we obtain $n$ distinct $(4n)$-cycles with the same mass
center.
\end{Rmk}

\begin{Rmk}
Let $x^{(m)}$ denote the $(4n)$-periodic point with itinerary
$\gamma^{(m)}$ for $0\leq m < n$. Then it can be calculated that
\begin{eqnarray*}
x^{(0)} &=& \frac{2}{3(2^{2n}+1)} + \frac{3\cdot
2^{4n-2}-2^{2n-2}}{2^{4n}-1};\\
 x^{(m)} &=& \frac{2}{3(2^{2n}+1)} +
\frac{2^{2n-2m}}{2^{2n}-1} \ \ \ \ \hbox{for $m=1,\cdots,n-1$.}
\end{eqnarray*}

Using equations (\ref{eqn:sym}), (\ref{eqn:cm}) and the common
trace vector of $x^{(m)}$'s, one can calculate the common mass
center of their orbits:
\begin{eqnarray*}
 \bar{x} & = & \frac{1}{4n (2^{4n}-1)}
               \left[2^{2n}s_{2n} +
                     \sum_{j=1\atop{j\;\textrm{is odd}}}^{4n-1}2^{4n-j}s_j +
                     \sum_{j=2\atop{j\;\textrm{is even}}}^{2n-2}2^{4n-j}s_j +
                     \sum_{j=2n+2\atop{j\;\textrm{is even}}}^{4n-2}2^{4n-j}s_j \right] \\
 & = &  \frac{1}{4n (2^{4n}-1)}
        \left[2n\left(2^{2n}+\sum_{r=0}^{2n-1}2^{2r+1}\right)
        +\sum_{r=1}^{n-1}\left((2r+4)2^{4n-2r}+(2n-2r+4)2^{2n-2r}\right)\right] \\
 & = & \frac{1}{3}+\frac{5}{9n}-\frac{13}{9n(2^{2n}+1)}-\frac{2}{n(2^{4n}-1)}.
\end{eqnarray*}

If we replace $n$ by $2n$ in these formulas, then we obtain the
$(8n)$-periodic points with itineraries $\beta^{(m)}$ and the
common mass center of their orbits.
\end{Rmk}

\section{Orbits without mass center and the set of mass centers}\label{sec:masscenters}

In this section we construct uncountably many orbits without mass center.
We also show that the set of all mass centers is the interval
$[0,\frac{2}{3}]$ and the set of mass centers of cycles is dense in $[0,\frac{2}{3}]$.
Moreover, for every $x$ in $[0,\frac{2}{3}]$, there are uncountably many orbits with mass
center $x$.

We start with some simple calculations. Given a positive integer
$k$, for any $0\leq i \leq k$ we have
\begin{eqnarray*}
T^i\left(\frac{2}{3}(\frac{1}{2^k}+\epsilon)\right) =
\frac{2}{3}\left(\frac{1}{2^{k-i}}+2^i\epsilon\right)\quad
\textrm{if}\;\;\frac{-1}{2^{k}} \leq \epsilon \leq \frac{1}{2^{k+1}},
\end{eqnarray*}
and
\begin{eqnarray*}
T^i\left(\frac{2}{3}(1+\epsilon)\right) =
\frac{2}{3}\left(1+(-1)^i 2^i\epsilon\right)\quad
\textrm{if}\;\;\frac{-1}{2^{k+1}} \leq \epsilon \leq \frac{1}{2^{k}}.
\end{eqnarray*}
In particular, both identities are valid for $|\epsilon| \leq \frac{1}{2^{k+1}}$.

Consider an infinite sequence of {\em even} integers $\{a_n\}_{n=0}^\infty$,
where $a_0=0$ and $a_n < a_{n+1}$ for all $n$. Let
\begin{eqnarray}\label{eqn:c}
c= \frac{2}{3}\sum_{k=0}^\infty \frac{(-1)^k}{2^{a_k}}.
\end{eqnarray}
Then
\begin{eqnarray*}
T^{1}(c) &=& \frac{2}{3}\left(1+\sum_{k\geq 1}\frac{(-1)^{k-1}}{2^{a_k-1}}\right) \\
  &\vdots &  \\
T^{a_1}(c) &=& \frac{2}{3}\left(1+\sum_{k\geq 1}\frac{(-1)^{k-a_1}}{2^{a_k-a_1}}\right)
         \;=\; \frac{2}{3}\sum_{k\geq 2}\frac{(-1)^k}{2^{a_k-a_1}} \\
T^{a_1+1}(c) &=& \frac{2}{3}\sum_{k\geq 2}\frac{(-1)^k}{2^{a_k-a_1-1}} \\
  &\vdots &  \\
T^{a_2}(c) &=& \frac{2}{3}\sum_{k\geq 2}\frac{(-1)^k}{2^{a_k-a_2}}
         \;=\; \frac{2}{3}\left(1+\sum_{k\geq 3}\frac{(-1)^{k-a_2}}{2^{a_k-a_2}}\right).
\end{eqnarray*}
The $a_k$'s are assumed to be even so the signs come out just right.
Proceed inductively, then for any $m\geq 0$ we have
\begin{eqnarray*}
T^n (c) &=& \frac{2}{3}
\left(1+\sum_{k\geq 2m+1}\frac{(-1)^{k-n}}{2^{a_k-n}}\right) \ \ \ \ \textrm{if}\
a_{2m}\leq n \leq a_{2m+1}
\end{eqnarray*}
and
\begin{eqnarray*}
T^n (c) &=& \frac{2}{3}\sum_{k\geq 2m+2}\frac{(-1)^k}{2^{a_k-n}} \
\ \ \ \textrm{if}\  a_{2m+1}\leq n \leq a_{2m+2}.
\end{eqnarray*}
Therefore,
\begin{eqnarray*}
\sum_{n=a_{2m}}^{a_{2m+1}-1}T^n (c)
&=& \frac{2}{3}\left[(a_{2m+1}-a_{2m})+
    \left(\sum_{n=a_{2m}}^{a_{2m+1}-1}\frac{(-1)^{n+1}}{2^{a_{2m+1}-n}} \right)
    \left(\sum_{k\geq 2m+1}\frac{(-1)^{k-1}}{2^{a_k-a_{2m+1}}}\right)\right] \\
&=& \frac{2}{3}\left[ (a_{2m+1}-a_{2m})+ \frac{1}{3}\left(1-\frac{1}{2^{a_{2m+1}-a_{2m}}}\right)
    \left(\sum_{k\geq 2m+1}\frac{(-1)^{k-1}}{2^{a_k-a_{2m+1}}}\right)\right] \\
\sum_{n=a_{2m+1}}^{a_{2m+2}-1}T^n (c)
&=& \frac{2}{3}\left(\sum_{n=a_{2m+1}}^{a_{2m+2}-1}\frac{1}{2^{a_{2m+2}-n}} \right)
    \left(\sum_{k\geq 2m+2}\frac{(-1)^k}{2^{a_k-a_{2m+2}}}\right) \\
&=& \frac{2}{3}\left(1- \frac{1}{2^{a_{2m+2}-a_{2m+1}}} \right)
    \left(\sum_{k\geq 2m+2}\frac{(-1)^k}{2^{a_k-a_{2m+2}}}\right).
\end{eqnarray*}
It follows that
\begin{eqnarray} \label{eqn:T-even}
\frac{2}{3} (a_{2m+1}-a_{2m}) < \sum_{n=a_{2m}}^{a_{2m+1}-1}T^n
(c) < \frac{2}{3} (a_{2m+1}-a_{2m})+1
\end{eqnarray}
and
\begin{eqnarray} \label{eqn:T-odd}
0< \sum_{n=a_{2m+1}}^{a_{2m+2}-1}T^n (c) < 1.
\end{eqnarray}

We can now make use of these calculations to prove Theorem~\ref{thm:nocenter}.
Let $\{a_n\}_{n=0}^\infty$ be a strictly increasing sequence of even integers with $a_0=0$ and
$a_{2k}=10^{2k}$, $a_{2k+1}=10^{2k+1}$ for infinitely many $k$'s.
For such a $k$,
\begin{eqnarray*}
\frac{1}{10^{2k}}\sum_{n=0}^{10^{2k}-1}T^n(c) <
\frac{1}{10^{2k}}\left({10^{2k-1}+\sum_{n=10^{2k-1}}^{10^{2k}-1}T^n(c)}\right)
< \frac{10^{2k-1} +1}{10^{2k}} \leq \frac{1}{5}.
\end{eqnarray*}
On the other hand
\begin{eqnarray*}
\frac{1}{10^{2k+1}}\sum_{n=0}^{10^{2k+1}-1}T^n(c) >
\frac{1}{10^{2k+1}}\sum_{n=10^{2k}}^{10^{2k+1}-1}T^n(c)
> \frac{\frac{2}{3}(10^{2k+1}-10^{2k})}{10^{2k+1}} = \frac{3}{5}.
\end{eqnarray*}
Thus the orbit of $c$ has no mass center.
It can be easily seen that there are uncountably many sequences $\{a_n\}_{n=0}^\infty$ with the
prescribed properties, where hence their corresponding points in (\ref{eqn:c})
generate uncountably many distinct orbits (since every orbit is countable).
This proves Theorem~\ref{thm:nocenter}.

If $\frac{2}{3}\leq x \leq 1$, then the average
$$
\frac{x+T(x)}{2} = 1-\frac{x}{2} \leq \frac{2}{3}.
$$
From this it follows easily that the mass center of an orbit can
be at most $\frac{2}{3}$, if it exists. We are now ready to show
that mass centers of cycles are dense in the interval $[0,
\frac{2}{3}]$. It suffices to show the following: for any rational
number $r\in[0,\frac{2}{3}]$ and any $\epsilon >0$, there exists a
periodic-point $c$ such that $|\bar{c}-r|<\epsilon$ where
$\bar{c}$ is the mass center of the periodic orbit of $c$. We may
assume that $r\neq \frac{2}{3}$. Let $r=\frac{p}{q}$ where $p$ and
$q$ are positive integers such that $3p< 2q$. Fix a positive {\em
even} integer $t$, for $m\geq 0$ let
\begin{eqnarray*}
a_{2m}=m(2qt)\ \ \text{and} \ \ a_{2m+1}=3pt + m(2qt).
\end{eqnarray*}
Then $\{a_n\}$ is a strictly increasing sequence of even integers
with $a_0 = 0$, and
\begin{eqnarray*}
c=\frac{2}{3}\sum_{i=0}^\infty \frac{(-1)^i}{2^{a_i}}
=\frac{2}{3}\left(1-\frac{1}{2^{a_1}}\right)
\sum_{k=0}^{\infty}\frac{1}{2^{k
a_2}}=\frac{2^{a_2-a_1+1}(2^{a_1}-1)}{3(2^{a_2}-1)}.
\end{eqnarray*}
Since $a_1=3pt$ is even, $(2^{a_1}-1)$ is divisable by $3$.
Therefore $c$ is an $a_2$-periodic point, and
\begin{eqnarray*}
\bar{c} = \frac{1}{a_2}\sum_{n=0}^{a_2-1} T^n(c).
\end{eqnarray*}
Now we use the inequalities (\ref{eqn:T-even}) and (\ref{eqn:T-odd}), then
\begin{eqnarray*}
\frac{2}{3} (a_1-0) < \sum_{n=0}^{a_1-1}T^n (c) < \frac{2}{3}
(a_1-0)+1
\end{eqnarray*}
and
\begin{eqnarray*}
0< \sum_{n=a_1}^{a_2-1}T^n (c) < 1.
\end{eqnarray*}
Combining them we obtain
\begin{eqnarray*}
\frac{p}{q}=\frac{\frac{2}{3}a_1}{a_2} <
\frac{1}{a_2}\sum_{n=0}^{a_2-1} T^n(c) <
\frac{(\frac{2}{3}a_1+1)+1}{a_2}=\frac{p}{q} + \frac{1}{qt}.
\end{eqnarray*}
Choose $t$ so that $\frac{1}{qt}< \epsilon$, then
\begin{eqnarray*}
|\bar{c}-r|=\left|\frac{1}{a_2}\sum_{n=0}^{a_2-1}
T^n(c)-\frac{p}{q}\right|< \epsilon.
\end{eqnarray*}

We now show that for every $x$ in the interval $[0, \frac{2}{3}]$,
there are uncountably many orbits with mass center $x$. First we
deal with the case of $x < \frac{2}{3}$. Let $x_n = \lceil nx
\rceil + e_n$, where $e_n = 1$ or 2; that is, $x_n$ is one of the
two possible integers such that
\begin{eqnarray*}
nx+1 \leq x_n < nx +3.
\end{eqnarray*}
Then there exists a positive integer $n_0$ such that $0< 3 x_n <
2n$ if $n\geq n_0$. For $m\geq 0$ let
\begin{eqnarray*}
a_{2m}=4\sum_{k=0}^{m-1} (k+n_0) = 2m(m+2n_0-1)\ \ \text{and} \ \
a_{2m+1}=a_{2m} + 6 x_{m+n_0}.
\end{eqnarray*}
Then $\{a_n\}$ is a strictly increasing sequence of even integers
with $a_0=0$, so we may define $c$ by (\ref{eqn:c}). Combining
inequalities (\ref{eqn:T-even}) and (\ref{eqn:T-odd}) we obtain
\begin{eqnarray*}
4 x_{m+n_0} \;<\; \sum_{n= a_{2m}}^{a_{2m+2}-1} T^n (c) \;<\; 4 x_{m+n_0}
+2.
\end{eqnarray*}
For any positive integer $N$, there exists $M$ such that $a_{2M}
\leq N < a_{2M+2}$. Thus we have
\begin{eqnarray*}
\sum_{n=0}^{N-1} T^n (c) &\geq& \sum_{m=0}^{M-1}
\sum_{n=a_{2m}}^{a_{2m+2}-1} T^n (c) \;>\; 4\sum_{m=0}^{M-1} x_{m+n_0}\\
 &\geq& 4\sum_{m=0}^{M-1}[ (m + n_0)x + 1] \;=\;
 4x\sum_{m=0}^{M-1} (m + n_0) + 4M \;=\; x a_{2M} + 4M
\end{eqnarray*}
and
\begin{eqnarray*}
\sum_{n=0}^{N-1} T^n (c) &\leq& \sum_{m=0}^{M}
\sum_{n=a_{2m}}^{a_{2m+2}-1} T^n (c) \;<\; 4\sum_{m=0}^{M} x_{m+n_0} +
2(M+1) \\
 &\leq& 4\sum_{m=0}^{M} [(m + n_0)x+3] + 2(M+1) \;=\;
 x a_{2M+2} +14(M+1).
\end{eqnarray*}
Hence
\begin{eqnarray*}
x \frac{a_{2M}}{N} + \frac{4M}{N} \;<\;
 \frac{1}{N}\sum_{n=0}^{N-1} T^n (c) \;<\;
 x \frac{a_{2M+2}}{N} +\frac{14(M+1)}{N}.
\end{eqnarray*}
Note that $a_{2M}$, $N$ and $a_{2M+2}$ are all of order $2M^2(1 +
o(1))$. Let $N \rightarrow \infty$ and we obtain $\bar{c} = x$ as
desired. Since we have two choices for each $x_n$, this produces
uncountably many different $c$.

It remains to deal with the case $x=\frac{2}{3}$. In this case let
$x_n = \lceil nx \rceil - e_n$, where $e_n = 1$ or 2. Once again
there exists a positive integer $n_0$ such that $0 < 3 x_n < 2n$
if $n \geq n_0$. The rest is pretty much the same as before. This
concludes the proof of Theorem~\ref{thm:dense}.

\section{Further examples}\label{sec:ex}

There are several patterns of itineraries
that can be suitably permuted without altering their trace vectors.
We will discuss one of them to which $\beta^{(1)}$ and $\gamma^{(1)}$
in section~\ref{sec:family} belong.

To simplify our notations, we will denote itineraries by row
vectors instead of column vectors in this section. For any string
$\mathbf{u}=[u_1,\cdots,u_k]$ consisting of $0$'s and $1$'s,
$\mathbf{u}^{-1}$ denotes the string obtained by reversing the
order of $u_i$'s in $\mathbf{u}$, and $|\mathbf{u}|$ denotes the
sum of $u_i$'s over $\Z/2\Z$. All itineraries are treated as
cyclic vectors.

\begin{Prop}\label{prop0011}
Let $\beta$ be an itinerary of the form
\begin{eqnarray*}
 \beta & = & [0,\p,0,\q,
              1,\p^{-1},1,\r]
\end{eqnarray*}
where $\p=[p_1,\cdots, p_j]$, $\q=[q_1,\cdots,q_k]$ and
$\r=[r_1,\cdots,r_k]$
are subject to conditions
\begin{eqnarray}\label{eqn:cond1}
 \q = \q^{-1},\quad \r = \r^{-1},
 \quad |\q| = |\r|.
\end{eqnarray}
Then the trace vector of $\beta$ does not change if $\p$ and
$\p^{-1}$ are switched; that is, $s(\beta)=s(\gamma)$, where
\begin{eqnarray*}
 \gamma & = & [0,\p^{-1},0,\q,1,\p,1,\r].
\end{eqnarray*}
\end{Prop}

Note that the presence of the two $0$'s and $1$'s implies that
$\gamma$ is in general not in the orbit of $\beta$ under the
action of the dihedral group (\ref{eqn:dihedral}). We will call
the two 0's and 1's the {\em four corners} of $\beta$ and
$\gamma$.
\begin{proof}
Let $n=2j+2k+4$ denote the length of $\beta$. By (\ref{eqn:sym}),
we only need to show that $s_l(\beta)=s_l(\gamma)$ for $1\leq l
\leq n/2$. Let us fix such an $l$. Define
$$
\begin{array}{ccc}
\beta'=[0,\p^{-1},0,\r,1,\p,1,\q], &
\beta''=[1,\p,1,\q,0,\p^{-1},0,\r], &
\beta'''=[1,\p^{-1},1,\r,0,\p,0,\q]; \vspace{1mm}\\
\gamma'=[0,\p,0,\r,1,\p^{-1},1,\q], &
\gamma''=[1,\p^{-1},1,\q,0,\p,0,\r], &
\gamma'''=[1,\p,1,\r,0,\p^{-1},0,\q].\\
\end{array}
$$
Then $\beta',\beta''$ and $\beta'''$ are in the same orbit of
$\beta$ under the action of the dihedral group
(\ref{eqn:dihedral}), so they all have the same trace vector.
Similarly $\gamma',\gamma''$ and $\gamma'''$ are in the same orbit
of $\gamma$, so they all have the same trace vector as well. For a
cyclic vector $\mathbf{x}=[x_1,\cdots,x_n]$, let
$\widehat{\mathbf{x}}=[x_i,\cdots,x_{i+l-1}]$ where it is
understood that the subindices are in $\Z/n\Z$. Since $l\leq n/2$,
the substring $\widehat{\beta}$ covers at most two corners of
$\beta$.

{\em Case 1: $\widehat{\beta}$ covers zero corners.} If
$\widehat{\beta}$ covers nothing of $\p$ and $\p^{-1}$, then note
that $\gamma$,$\gamma'$,$\gamma''$ and $\gamma'''$ are obtained
from $\beta$,$\beta'$,$\beta''$ and $\beta'''$ by switching $\p$
and $\p^{-1}$. If $\widehat{\beta}$ covers nothing of $\q$ and
$\r$, then note that $\gamma$,$\gamma'$,$\gamma''$ and $\gamma'''$
are obtained from $\beta'$,$\beta$,$\beta'''$ and $\beta''$ by
switching $\q$ and $\r$. In either case we have
$\{|\widehat{\beta}|,|\widehat{\beta'}|,|\widehat{\beta''}|,|\widehat{\beta'''}|\}$
and
$\{|\widehat{\gamma}|,|\widehat{\gamma'}|,|\widehat{\gamma''}|,|\widehat{\gamma'''}|\}$
are the same as multisets.

{\em Case 2: $\widehat{\beta}$ covers one corner.} Notice that if
we switch 0 and 1 in the four corners, then $\beta$ becomes
$\beta''$ and $\beta'$ becomes $\beta'''$. Therefore as multisets
$\{|\widehat{\beta}|,|\widehat{\beta'}|,|\widehat{\beta''}|,|\widehat{\beta'''}|\}
= \{0,0,1,1\}$. For similar reason we also have
$\{|\widehat{\gamma}|,|\widehat{\gamma'}|,|\widehat{\gamma''}|,|\widehat{\gamma'''}|\}
= \{0,0,1,1\}$.

{\em Case 3: $\widehat{\beta}$ covers two corners.} The argument
here is similar to Case 1. If $\widehat{\beta}$ covers a complete
$\p$ (or $\p^{-1}$), then it covers nothing of $\p^{-1}$ (or $\p$)
since $l\leq n/2$. In this case switch $\p$ and $\p^{-1}$. If
$\widehat{\beta}$ covers a complete $\q$ (or $\r$), then it covers
nothing of $\r$ (or $\q$). In this case switch $\q$ and $\r$.
Since $|\p^{-1}|=|\p|$ and $|\q|=|\r|$, in either case we have
$\{|\widehat{\beta}|,|\widehat{\beta'}|,|\widehat{\beta''}|,|\widehat{\beta'''}|\}=\{|\widehat{\gamma}|,|\widehat{\gamma'}|,|\widehat{\gamma''}|,|\widehat{\gamma'''}|\}$.

Now let $i$ run through $\{1, \cdots, n\}$, then we have
$4s_l(\beta)=4s_l(\gamma)$.
\end{proof}

\begin{Rmk}
Same arguments can be applied to itineraries of the form
\begin{eqnarray*}
 \beta & = & [0,\p,1,\q,
              0,\p^{-1},1,\r]
\end{eqnarray*}
where $\p=[p_1,\cdots, p_j]$, $\q=[q_1,\cdots,q_k]$ and
$\r=[r_1,\cdots,r_k]$
are subject to conditions (\ref{eqn:cond1}).
\end{Rmk}

\begin{Exam}
Both the $\beta^{(1)}$ in (\ref{eqn:betam}) and $\gamma^{(1)}$ in (\ref{eqn:gammam})
have the structure described in Proposition~\ref{prop0011}.
For instance, the $\beta^{(1)}$ in (\ref{eqn:beta123}) as a cyclic vector can be written
\begin{eqnarray*}
 \beta^{(1)} & = & [0,\underbrace{1,0}_{\p},0,\underbrace{1,1,1,1,1,1,1,1}_{\q},
                    1,\underbrace{0,1}_{\p^{-1}},1,\underbrace{0,0,0,0,0,0,0,0}_{\r}].
\end{eqnarray*}
By Proposition~\ref{prop0011}, the itinerary
\begin{eqnarray*}
 \beta^{(0)} & = & [0,\underbrace{0,1}_{\p^{-1}},0,\underbrace{1,1,1,1,1,1,1,1}_{\q},
                    1,\underbrace{1,0}_{\p},1,\underbrace{0,0,0,0,0,0,0,0}_{\r}]
\end{eqnarray*}
has the same trace vector as $\beta^{(1)}$. Their trace vector is
\begin{eqnarray*}
 (12,6,12,8,12,10,12,12,12,14,12,12,12,14,12,12,12,10,12,8,12,6,12,0).
\end{eqnarray*}
Proposition~\ref{prop0011} also gives another reason why the
$\gamma^{(0)}$ and $\gamma^{(1)}$ in (\ref{eqn:gammam}) have the
same trace vector.
\end{Exam}

\begin{Exam}
Ideally, one would hope to have itineraries where
Proposition~\ref{prop0011} could be applied multiple times to
produce many distinct cycles. Here we present such an interesting
example for the case $n=20$. Unfortunately it is unclear how to
generalize it to general $n$. Let
\begin{eqnarray*}
 \beta^{(0)} & = & [0,\underbrace{0,0,0,1,1,0}_{\p},0,\underbrace{0,0}_{\q},
                    1,\underbrace{0,1,1,0,0,0}_{\p^{-1}},1,\underbrace{1,1}_{\r}]
\end{eqnarray*}
After reversing $p_i$'s and shifting the cyclic vector, the new itinerary
$\beta^{(1)}$ also has the pattern studied in Proposition~\ref{prop0011}.
By repeating this process, we obtain
\begin{eqnarray*}
 \beta^{(1)} & = & [0,\underbrace{1,0,0}_{\p},0,\underbrace{1,1,0,1,1}_{\q},
                    1,\underbrace{0,0,1}_{\p^{-1}},1,\underbrace{0,0,0,0,0}_{\r}] \\
 \beta^{(2)} & = & [0,\underbrace{1,0}_{\p},0,\underbrace{0,0,0,0,0,0}_{\q},
                    1,\underbrace{0,1}_{\p^{-1}},1,\underbrace{0,1,1,1,1,0}_{\r}] \\
 \beta^{(3)} & = & [0,\underbrace{0,0,0,1,1}_{\p},0,\underbrace{1,0,1}_{\q},
                    1,\underbrace{1,1,0,0,0}_{\p^{-1}},1,\underbrace{0,0,0}_{\r}].
\end{eqnarray*}
We resume $\beta^{(0)}$ by reversing the $p_i$'s in $\beta^{(3)}$.
It is easy to see that these itineraries and their mirror itineraries are all different, yielding eight
distinct cycles with the same mass center.
In the order of $\beta^{(0)},\cdots,\beta^{(3)}$ followed by their mirror itineraries,
they correspond respectively to
$$
\frac{33658}{1048575},\;\frac{504896}{1048575},\;\frac{523412}{1048575},\;\frac{39664}{1048575},\;
\frac{776224}{1048575},\;\frac{17852}{1048575},\;\frac{337916}{1048575},\;\frac{125728}{1048575}.
$$
Their common trace vector is
\begin{eqnarray*}
 (8,8,8,10,10,10,12,12,12,12,12,12,12,10,10,10,8,8,8,0),
\end{eqnarray*}
and their common center of mass is $\frac{2170804}{5242875}\approx0.414048$.
\end{Exam}

\begin{Exam}\label{ex-24}
There are patterns of itineraries that allow much more distinct
cycles with the same mass centers than the one constructed in
section~\ref{sec:family}. However it is unclear how to generalize these
patterns. The following itineraries along with their mirror
itineraries correspond $16$ cycles of period $24$ with the same
center of mass.
\begin{eqnarray*}
 \beta^{(0)} & = & [0,0,0,0,0,0,1,1,1,0,1,0,1,0,0,1,1,1,0,1,1,1,0,1]\\
 \beta^{(1)} & = & [0,0,0,0,0,1,1,0,1,1,0,1,1,1,1,1,0,1,0,1,0,0,0,1]\\
 \beta^{(2)} & = & [0,0,0,0,0,1,1,0,1,1,1,0,1,1,0,0,1,0,0,1,1,1,0,1]\\
 \beta^{(3)} & = & [0,0,0,0,0,1,1,1,1,0,1,0,0,0,1,1,1,0,1,1,0,1,0,1]\\
 \beta^{(4)} & = & [0,0,0,0,0,1,0,1,0,0,0,1,1,1,1,1,0,1,1,0,1,1,0,1]\\
 \beta^{(5)} & = & [0,0,0,0,1,1,0,0,0,1,1,1,0,1,1,1,0,1,0,1,1,0,0,1]\\
 \beta^{(6)} & = & [0,0,0,0,1,1,0,1,1,0,0,1,0,1,1,1,0,1,0,0,0,1,1,1]\\
 \beta^{(7)} & = & [0,0,0,0,1,0,1,1,0,1,1,0,0,0,1,1,0,1,1,1,1,0,0,1].
\end{eqnarray*}
From $\beta^{(0)}$ to $\beta^{(7)}$, they are respectively itineraries of
$$
\frac{183958}{16777215},\;\frac{300446}{16777215},\;\frac{309014}{16777215},\;\frac{343334}{16777215},\;
\frac{398774}{16777215},\;\frac{547438}{16777215},\;\frac{596602}{16777215},\;\frac{900526}{16777215}.
$$
The mirror itineraries of $\beta^{(0)},\cdots,\beta^{(7)}$ correspond respectively to
$$
\frac{13821568}{16777215},\;\frac{15946304}{16777215},\;\frac{13752896}{16777215},\;\frac{13203776}{16777215},\;
\frac{14373056}{16777215},\;\frac{15512608}{16777215},\;\frac{12366112}{16777215},\;\frac{15432544}{16777215}.
$$
Their common trace vector is
\begin{eqnarray*}
 (12,12,10,14,12,12,12,10,14,12,12,16,12,12,14,10,12,12,12,14,10,12,12,0),
\end{eqnarray*}
and their common center of mass is $\frac{24897599}{50331645}\approx0.494671$.
\end{Exam}

\begin{Exam}\label{counterexample}
The criterion in section~\ref{sec:binary} for two distinct cycles to have the
same mass center is sufficient but not necessary.
For example, the fractions $\frac{5414}{131071}$ and $\frac{10090}{131071}$
are points of period $17$ with itineraries
\begin{eqnarray*}
 \beta  & = & [0,0,0,0,1,1,1,1,1,1,0,1,1,0,1,0,1]\\
 \gamma & = & [0,0,0,1,1,0,1,0,0,1,1,0,1,1,1,1,1].
\end{eqnarray*}
Their trace vectors are respectively
\begin{eqnarray*}
 s(\beta)  & = & (10,8,8,8,10,8,6,8,8,6,8,10,8,8,8,10,0)\\
 s(\gamma) & = & (10,8,8,8,6,14,8,10,10,8,14,6,8,8,8,10,0).
\end{eqnarray*}
One can easily check that $\frac{1185588}{2228207}\approx0.532082$ is their common center of mass.
\end{Exam}

\begin{Exam}\label{ex:last}
Cycles of different lengths may have the same center of mass.
For instance, the orbit of $\frac{46}{275}$ has length $20$ and the orbit of
$\frac{1158}{7735}$ has length $24$.
One can easily verify that their mass centers are both equal to $\frac{3}{5}$.

We have mentioned in section~\ref{sec:intr} that, by the Birkhoff
ergodic theorem, the mass center of generic orbits is
$\frac{1}{2}$, the space average or integral of the tent map.
Periodic orbits are exceptional but there are actually examples of
cycles where the mass center is also $\frac{1}{2}$. For instance,
$$\frac{19}{18564},\;\frac{17}{14460},\;\frac{13}{4820},\;\frac{151}{55692}$$
are all $24$-periodic points and $\frac{1}{2}$ is the common mass center of their orbits.
\end{Exam}

\section{Final Remarks and Related Problems}\label{sec:rmk}

The discussions presented so far are focused on periodic points of
the form $\frac{2m}{2^n-1}$. They are also applicable to periodic
points of the form $\frac{2m}{2^n+1}$ because
$\frac{2m}{2^n+1}=\frac{2m(2^n-1)}{2^{2n}-1}$. We say an $n$-cycle
is of the {\em first type} if it is generated by points of the
form $\frac{2m}{2^n-1}$, and it is of the {\em second type} if
otherwise. An $n$-periodic point of the second type is a
$2n$-periodic point of the first type. Two of the examples in the
introduction are of the second type.

Let $c(n)$ denote the maximum number of distinct first type $n$-cycles with
the same center of mass, and $d(n)$ denote the same number for second type $n$-cycles.
Exhaustive search by computer shows that $c(n)$ and $d(n)$ are
{\scriptsize
\begin{center}
\begin{tabular}{c|c|c|c|c|c|c|c|c|c|c|c|c|c|c|c|c|c|c|c|c|c|c}
$n$    & $1\sim6$ & $7\sim11$& $12$ & $13$& $14$ & $15$ & $16$ & $17$ & $18$ & $19$ & $20$
       & $21$ & $22$ & $23$  & $24$ & $25$& $26$ & $27$ & $28$ & $29$ & $30$ & $31$  \\ \hline
$c(n)$ & $1$      & $2$      & $4$  & $4$ & $4$  & $4$  & $6$  & $6$  & $6$  & $8$  & $8$
       & $8$  & $8$  & $8$   & $16$ & $10$& $10$ & $14$ & $38$ & $16$ & $24$ & $10$  \\ \hline
$d(n)$ & $1$      & $2$      & $4$  & $4$ & $5$  & $4$  & $4$  & $4$  & $8$  & $4$  & $8$
       & $8$  & $6$  & $8$   & $10$ & $8$ & $10$ & $14$ & $12$ & $16$ & $24$ & $10$
\end{tabular}
\end{center}
}

Let $e(n)$ be the maximum number of distinct $n$-cycles of either
type with the same center of mass. It turns out that $e(n) =
\max\{ c(n), d(n) \}$ for $n\leq 31$. Also note that in our
notation $n$-cycles include all cycles whose lengths are divisors
of $n$. If we only consider cycles whose lengths are exactly $n$,
then the values of $c(n)$, $d(n)$ and $e(n)$ remain the same for
$n \leq 31$. Are these observations true for arbitrary $n$?

Another interesting question about cycles with the same mass
center is how fast $e(n)$ grows. From section~\ref{sec:family} we
know that $e(4n) \geq n$. Is $e(n)$ of order $O(n)$ as $n$ goes to
infinity?
If so, what are the values of
$\liminf_{n\rightarrow\infty}\frac{e(n)}{n}$ and
$\limsup_{n\rightarrow\infty}\frac{e(n)}{n}$? Intuitively $e(n)$
grows faster on composite numbers because itineraries of
$n$-cycles with prime $n$ seems to admit fewer possible
trace-invariant permutations.
Note that the values of $e(29)$ and $e(31)$ are respectively less than half of
$e(28)$ and $e(30)$.
Each of $e(17)$ and $e(19)$ is
attained at only one example where not all of the distinct cycles
have the same trace vector.
How fast does $e(n)$ grow on prime numbers?

Let $C_n$ be the set of reals with the property that each point of
$C_n$ is the mass center of at least $n$ distinct cycles. It would
be interesting to understand the intersection of the nested
sequence $\{C_n\}$. Is this intersection nonempty? Equivalently,
is there a point which is the mass center of infinitely many
cycles?

Our work and the problems raised above are related to the class of exceptional orbits
to which the Birkhoff ergodic theorem does not apply.
Exceptional orbits are rare but their dynamics can be used to characterize different chaotic maps.
We have seen that, for the standard tent map,
the dynamics exhibited by the exceptional orbits are rich and maybe surprising.
Similar dynamical phenomena may occur for many other chaotic interval maps or even higher dimensional maps.

\hfill\newline
\noindent{\bf Acknowledgement.}\\
This research work is partially supported by the National Science
Council and the National Center for Theoretical Sciences in
Taiwan. Part of the work was completed while the second author was
a visiting researcher at the Center in 2007. The second author
would like to thank the Center for the support.
\smallskip


\end{document}